\documentclass[review, 12pt]{elsarticle}

\usepackage{graphicx}
\usepackage{amssymb}
\usepackage{amsthm}
\usepackage{amsmath}

\newtheorem{thm}{Theorem}
\newtheorem{alg}{Algorithm}
\newtheorem{ex}{Example}

\usepackage{lineno}

\journal{Chaos, Solitons \& Fractals}

\begin{document}

\begin{frontmatter}



\title{Testing of fractional Brownian motion in a noisy environment}



\author[ad]{Micha\l{} Balcerek\corref{cor1}}
\cortext[cor1]{Corresponding author}
\ead{michal.balcerek@pwr.edu.pl}

\author[ad]{Krzysztof Burnecki}

\address[ad]{Faculty of Pure and Applied Mathematics, Hugo Steinhaus Center, Wroc\l{}aw University of Science and Technology, Wyspia\'nskiego 27, 50-370 Wroc\l{}aw, Poland}

\begin{abstract}
Fractional Brownian motion (FBM) is related to the notions of self-similarity, ergodicity and long memory. 
These properties have made FBM important in modelling real-world phenomena in different experiments ranging from telecommunication to biology. However, these experiments are often disturbed by a noise which source can be, e.g., the instrument error. In this paper we propose a rigorous statistical test for FBM with added white Gaussian noise which is based on the autocovariance function. To this end we derive a distribution of the test statistic which is given explicitly by the generalized chi-squared distribution. This allows us to find critical regions for the test with a given significance level. We check the quality of the introduced test by studying its power for alternatives being FBM's with different self-similarity parameters and the scaled Brownian motion which is also Gaussian and self-similar. We note that the introduced test can be adapted to an arbitrary Gaussian process with a given covariance structure. 
\end{abstract}



\end{frontmatter}


\section{Introduction}
\label{S:1}
Fractional Brownian motion (FBM), introduced by Kolmogorov in 1940 \cite{kolmogorov1940wienersche,mandelbrot1968fractional}, is a generalisation of the classical Brownian motion (BM). Most of its statistical properties are characterized by the self-similarity parameter (Hurst exponent) $0<H<1$.
FBM, denoted by $B_H(t)$, is $H$-self-similar, namely for every $c>0$ we have $B_H(ct)\stackrel{D}{=}c^H B_H(t)$ in the sense of all finite dimensional distributions, and has stationary increments. It is the only Gaussian process satisfying these properties. With probability $1$, the graph of $B_H(t)$ has both Hausdorff dimension and box dimension of $2-H$.

The increments of FBM $Y_j=B_H(j+1)-B_H(j);\ j=0,1,\ldots$ are called fractional Gaussian noise (FGN). FGN has some remarkable properties. If $H=1/2$, then its autocovariance function (ACVF) $r(k)=0$ for $k\neq 0$ and hence it is the sequence
of independent random variables. The
situation is quite different when $H\neq 1/2$, namely the $Y_j$'s are correlated (dependent) and
the time series has the ACVF $r(k)$ of the power-law form:
\begin{equation}
\label{finbez3} 
r(k)\sim {\rm Var} Y_1\,
H(2H-1)k^{2H-2},\;\;\;\;{\rm as}\;\;k\rightarrow\infty.
\end{equation}
The $r(k)$ tends to $0$ as $k\rightarrow\infty$  for all
$H\neq 1/2$ much slower than exponentially. For example, the Ornstein-Uhlenbeck process, which is the  most common stationary process, has 
exponentially fast decaying correlations. Moreover, when $1/2<H<1$ $r(k)$ tends to zero so slowly that the sum
$\sum_{k=1}^{\infty}|r(k)|$ diverges to infinity. We say that in this case the increment process exhibits long memory (long-range dependence, persistence) \cite{beran2016long}. We also note that
the coefficient $H(2H-1)$ is positive, so the $r(j)$'s are positive for all large $j$, a behaviour referred to as ``positive dependence''.
Furthermore,
formula (\ref{finbez3}) by the Wiener Tauberian theorem \cite{zyg59} implies that the spectral density $h(\lambda)$ of has a pole at zero which leads to  a phenomenon often referred to as ``$1/f$ noise''. Such a behavior of the ACVF has become especially important in areas such as communication networks 
\cite{nor95, wiletal97} and finance \cite{bai96, Cheridito2003,FALLAHGOUL201723}.

If $0<H<1/2$, then $\sum_{k=1}^{\infty}|r(k)|<\infty$ and the spectral density tends
to zero as $|\lambda|\rightarrow 0$. Furthermore, the coefficient $H(2H-1)$ is negative and the $r(j)$'s
are negative for all large $j$. We say in that case that the sequence displays
negative power-like dependence called antipersistence, short or medium memory \cite{samorodnitsky2007long}. Antipersistence has been observed in financial time series for electricity price processes \cite{werprz00,wer06}, in climatology \cite{caretal07} and is widely pronounced in nanoscale biophysics in the context of viscoelastic systems \cite{kou08, magetal09, buretal12c, metzler2014anomalous, granik2019single}.


Since $\mathrm{E} B_H^2(t)=t^{2H}$, for $H\neq 1/2$ the second moment is not linear but sub- or super-linear.    
In physics, this behaviour is closely related to the notion of anomalous diffusion \cite{metzler2000random}, and
the second moment is called the (ensemble) mean-squared displacement (MSD). 
The sublinear form of MSD is related to subdiffusion  which is often observed in crowded systems, for example protein diffusion within cells, or diffusion through porous media, and the superlinear to superdiffusion \cite{ metzler2014anomalous}.

Since FBM has played an important role in many scientific disciplines and applied fields its proper identification and validation is an important issue.
In the literature different methods of estimating the self-similarity index $H$ have been developed \cite{taqtev98,douetal03,beran2016long}. 
An estimator based on MSD was studied in Ref. \cite{burgajsik11}.
It is an unbiased estimator with a very low variance and works remarkably well for the FBM \cite{hab}.
Rigorous statistical tests for the FBM already appeared in the literature, see \cite{sikora2017mean,sikora2018statistical}. They were based on the MSD and detrended moving average (DMA) statistics, respectively. 

Another stochastic process similar to FBM which recently attracted attention of physicists and mathematicians is the scaled Brownian motion (SBM) \cite{lim2002,bedrova2019,magdziarz2020}. The SBM $B_{s}(t)$ is a generalisation of the BM, namely $B_{s}(t)=B(t^{\alpha})$, $\alpha>0$, where $B(t)$ is the BM.  The process is Gaussian and self-similar like FBM (the Hurst exponent equals $\alpha/2$) but in general has independent and non-stationary increments so the memory structure of the increment process is completely different than that of FGN.

The FBM is observed in experiments which are often disturbed by a noise which source can be the  measurement and instrumentation error. We consider here a Gaussian white noise added to the FBM. This model was already studied, e.g., in Refs. \cite{burnecki2015estimating} where an idea of estimating both the self-similarity parameter and the magnitude of the Gaussian noise were presented. To the best of authors' knowledge there have been no statistical tests on such model in the literature.

The FGN is a moving average process, hence it is ergodic. 
Ergodicity is a very important characteristic since it is related to the Boltzmann hypothesis of equality of averages. It is also the reason, for which one often checks if the ensemble and time-average MSD's are coinciding. In case they are not, we say that the model (or system) exhibits weak ergodicity breaking \cite{bouchaud1992weak, jeon2011vivo}. Yet, equality of these values does not guarantee the ergodicity. Thus, in the literature several methods of checking the ergodicity have been developed \cite{magdziarz2011anomalous, janczura2015ergodicity, schwarzl2017quantifying,  slkezak2019codifference}. 
Convergence of the ACVF of a Gaussian process implies mixing which is a stronger property than ergodicity. Hence, in practice, it is often easier to check mixing since the ACVF of a Gaussian process is usually known.

The increment process of the FBM with added white Gaussian noise is also stationary and ergodic. Studying the ACVF of the process has been our motivation for developing a rigorous statistical test for the FBM with additive noise. We  believe that the ACVF statistic is simpler than statistics proposed for testing the pure FBM in Refs. \cite{sikora2017mean,sikora2018statistical}


In Section \ref{section:ergodicitydefinitions} we present the motivation for choosing the ACVF for testing of the FBM with noise, namely its simple structure and relation to ergodicity. We derive a distribution of the quadratic form corresponding to the ACVF for general Gaussian processes and specialize the result to the FBM with additive noise. In Section \ref{sec:test} we introduce a statistical test on FBM with additive noise based on ACVF. It can be used to test whether the data can be described by FBM with noise with given self-similarity and noise magnitude parameters. We show how to calculate quantiles of the statistic as a function of the parameters, which and we call a critical surface.
Section~\ref{sec:power} is devoted to the analysis of the power of the test. For alternatives we take FBM's and SBM's with varying self-similarity parameters.   We show the test can distinguish between models with a high efficiency. Section \ref{sec:con} concludes the article.

\section{Autocovariance function and its relation to ergodicity}
\label{section:ergodicitydefinitions}
We will recall here the relation between the ACVF and ergodicity of a stationary process. We will construct a statistical test for an arbitrary stationary Gaussian process based on the ACVF statistic.

Ergodicity is a very important property, because if stationary process $\left\{Y(t)\right\}_{t\geq0}$ is ergodic then the Boltzmann ergodic hypothesis is true, i.e. the time average converges to the ensemble average, that is
\begin{align}
\lim_{T\to\infty} \frac{1}{T} \int_0^T g(Y(t)) \mathrm{d}{t} = \mathbb{E} g(Y(0)),
\label{eq:ergdef}
\end{align}
{for any function $g$}, provided $\mathbb{E}|g(Y(0))|<\infty$ \cite{birkhoff1931proof, boltzmann2012theoretical}. Integration in (\ref{eq:ergdef}) is considered in the sense of trajectory by trajectory. This means that by observing one single trajectory we can infer characteristics of the whole system. For example, if $g(x) = x^2$, then the time averaged mean-squared displacement (TAMSD) of the process $\left\{Y(t)\right\}_{t\geq0}$ provides full information on the second moment of the process.

From Ref. \cite{maruyama1970infinitely} we know that a zero mean stationary Gaussian process $\left\{X(n)\right\}_{n = 0, 1, \ldots}$ with the ACVF $r(k)$ is mixing if and only if
  \begin{equation*}
  \lim_{k\to\infty} r(k) = 0. 
  \end{equation*}
    Mixing is a very important characteristic of a process and it it stronger property than ergodicity, namely the mixing process is also ergodic. Also, often it is much easier to check mixing than ergodicity. 
    
  The condition for ergodicity of Gaussian processes is also given in Ref. \cite{maruyama1970infinitely} where it states that a zero mean stationary Gaussian process $\left\{X(n)\right\}_{n = 0, 1, \ldots}$ with the ACVF $r(k)$ is ergodic if and only if
  \begin{equation*}
  \lim_{n\to\infty}\frac{1}{n} \sum_{k=1}^{n} |r(k)| = 0. 
  \end{equation*}

Therefore, to check either mixing or ergodicity of a zero mean stationary Gaussian process we should know a behavior of its ACVF.
In order to do so, we analyse the distribution of the sample ACVF estimator $\hat{r}(k)$:
  \begin{align}
      \hat{r}(k) = \frac{1}{N-k} \sum_{i=0}^{N-k-1} X_i X_{i+k}.
  \label{eq:rk}
  \end{align}
  The estimator can be written in a quadratic form:
  \begin{align*}
      \hat{r}(k) = \mathbf{X}^T A_k \mathbf{X}, \quad k=0, \ldots, N-1,
  \end{align*}
  where $\mathbf{X} = [X_0, X_1, \ldots, X_{N-1}]^T$, and matrix $A_k = \left[a_{i,j}\right]_{i=1,\ldots, N, j = 1, \ldots, N}$ is given by:
  \begin{align}
  a_{i,j} &= \begin{cases}
              \frac{1}{N-k} \frac{1}{2} & \textrm{ for every such pair that } |i-j|=k,\\
              0 & \textrm{otherwise,} 
              \end{cases}
              \label{eq:Amatrix}
  \end{align}
  if $k=1,2, \ldots, N-1$, and
  \begin{align*}
  a_{i,j} & = \begin{cases}
              \frac{1}{N} \quad \textrm{ if } i=j, \\
              0 \quad \textrm{otherwise.} 
              \end{cases}
  \end{align*}
  The matrix $A_k$'s has  non-zero elements on the diagonals starting in the $k$-th column and/or $k$-th row. For a better understanding of its construction, in the following example we present the matrices $A_1$ and $A_2$.
  \begin{ex}
Let us consider the matrices $A_k$ for $k=1$ and $k=2$.
      \begin{itemize}
          \item Case $k=1$. Matrix $A_1$ has the following form:
              \begin{align*}
                  (N-1) A_1 = 
                  \begin{bmatrix}
                      0 & \frac{1}{2} & 0 & 0 & 0 & 0 & \ldots & 0 \\
                      \frac{1}{2} & 0 & \frac{1}{2} & 0 & 0 & 0 & \ldots & 0 \\
                      0 & \frac{1}{2} & 0 & \frac{1}{2} & 0 & 0 & \ldots & 0\\
                      0 & 0 & \frac{1}{2} & 0 & \frac{1}{2} & 0 & \ldots & 0\\
                     \vdots & & & & & & \ddots &  \\
                      0 & 0 & \ldots & 0 & 0 & 0 & \frac{1}{2} & 0
                  \end{bmatrix}
              \end{align*}
          \item Case $k=2$. Matrix $A_2$ has the following form:
              \begin{align*}
                  (N-2) A_2 = 
                  \begin{bmatrix}
                      0 & 0 & \frac{1}{2} & 0 & 0 & 0 & \ldots & 0 \\
                      0 & 0 & 0 & \frac{1}{2} & 0 & 0 & \ldots & 0 \\
                      \frac{1}{2} & 0 & 0 & 0 & \frac{1}{2} & 0 & \ldots & 0\\
                      0 & \frac{1}{2}& 0 & 0 & 0 & \frac{1}{2} & \ldots & 0\\
                     \vdots & & & & & & \ddots &  \\
                      0 & 0 & \ldots & 0 & 0 & \frac{1}{2} & 0 & 0
                  \end{bmatrix}
              \end{align*}
      \end{itemize}
\end{ex}

 We now state the main result of this paper. 
 \begin{thm}
 \label{th:Q1}
 The quadratic form $Q(\mathbf{X})$ corresponding to the sample ACVF of the vector $\mathbf{X}$, namely
      \begin{align}
          Q(\mathbf{X}) \stackrel{df}{=} \frac{1}{N-k} \sum_{i=0}^{N-1-k} X_i X_{i+k} = \mathbf{X}^T A_k \mathbf{X}, \quad k=0, \ldots, N-1,
          \label{eq:Q1def}
      \end{align}
  has a generalized $\chi^2$ distribution, i.e. 
      \begin{align}
          Q(\mathbf{X}) \stackrel{D}{=} \sum_{j=1}^N \lambda^{(k)}_j U_j^2,
          \label{eq:Q1}
      \end{align}
  where $U_j^2$'s are independent random variables having the $\chi^2$ distribution with one degree of freedom, values $\left\{ \lambda^{(k)}_j \right\}_{j=1,\ldots,N}$ are eigenvalues of the matrix $\Sigma^{1/2} A_k \Sigma^{1/2}$, where the matrix $A_k$ is defined by formula (\ref{eq:Amatrix}), and the matrix $\Sigma$ is the covariance matrix of $\mathbf{X}$.
\begin{proof}
  Let us introduce the notation $\mathbf{Y} = \Sigma^{-1/2} \mathbf{X} \sim \mathbf{\mathcal{N}}(0, \mathbf{I}_N)$. Then, 
  \begin{align*}
      \mathbf{X}^T A \mathbf{X} = \mathbf{Y}^T \Sigma^{1/2} A \Sigma^{1/2} \mathbf{Y}.
  \end{align*}
  Based on the matrix spectral theorem \cite{meyer2000matrix}, we find the decomposition:
  \begin{align*}
      \Sigma^{1/2} A_k \Sigma^{1/2} = P^T \Lambda_k P, 
  \end{align*}
  where $P$ is an orthogonal matrix $P^T P = P P^T = \mathbf{I}$ and $\Lambda_k$ is a diagonal matrix with elements $\left\{ \lambda^{(k)}_j \right\}$ on the main diagonal. Those elements are eigenvalues of $\Sigma^{1/2} A_k \Sigma^{1/2}$.
  
	For $U = P \mathbf{Y} \sim \mathbf{\mathcal{N}}(0, \mathbf{I}_N)$ we have
  \begin{align*}
      Q(\mathbf{X}) & = \hat{r}(k) = \mathbf{X}^T A_k \mathbf{X} = \mathbf{Y}^T \Sigma^{1/2} A_k \Sigma^{1/2} \mathbf{Y} =\\
       & = \mathbf{Y}^T P_k^T \Lambda_k P_k \mathbf{Y} = (P_k \mathbf{Y})^T \Lambda_k (P_k \mathbf{Y}) =\\
       & = \mathbf{U}^T \Lambda_k \mathbf{U} = \sum_{j=1}^N \lambda^{(k)}_j U_j^2,
  \end{align*}
  where $U_j^2$'s are independent random variables having the $\chi^2$ distribution with one degree of freedom.
\end{proof}  
\end{thm}




\subsection{Application to the FBM with  noise}
In practice, the observed data are often disturbed by a noise (e.g. a measurement noise). We consider here the process given by:
  \begin{equation}
  X_H(t) = B_H(t) + \xi(t), 
	\label{eq:toyModelContinuousTime}
  \end{equation}
  where $\left\{ B_H(t) \right\}_{t\geq0}$ is the FBM and $\left\{ \xi(t) \right\}_{t\geq0}$ is the white Gaussian noise with variance $\sigma^2$.
  
  The increment process  $M(n) = X_H(n+1)- X_H(n), n = 0, 1, \ldots$ is a zero-mean stationary Gaussian process with the ACVF:
  \begin{equation}
    r_M(k) = 
    \begin{cases}
    r(k) + 2 \sigma^2 & \textrm{if } k=0,\\
    r(k) - \sigma^2 & \textrm{if } k=1,\\
    r(k) & \textrm{if } k>1,\\
    \end{cases}
    \label{eq:toymodelcov}
\end{equation}
where $r(k)$ is the ACVF of the increments of the FBM, i.e. 
\begin{equation}
r(k) = \frac{1}{2}\left((k+1)^{2H} + |k-1|^{2H} - 2 k^{2H} \right), \quad k = 0, 1, \ldots
\label{eq:fbmacf}
\end{equation}

We apply now Theorem \ref{th:Q1} to the FBM with added noise.
  For the considered process, the autocovariance matrix for the vector of its increment process has the following form:
  \begin{align}
          \Sigma =\begin{bmatrix}
                      r_M(0) & r_M(1) & r_M(2) & \ldots & r_M(N-1) \\
                      r_M(1) & r_M(0) & r_M(1) & \ldots & r_M(N-2) \\
                      r_M(2) & r_M(1) & r_M(0) & \ldots & r_M(N-3) \\
                      \vdots &        &        & \ddots & \vdots \\
                      r_M(N-1) & r_M(N-2) & r_M(N-3) & \ldots & r_M(0) \\
                  \end{bmatrix},
      \label{eq:covMatrix}
  \end{align}
  that is $\Sigma = [a_{ij}]_{i,j}$, where $a_{ij} = r_M(|i-j|)$, and $r_M(\cdot)$ is the ACVF of the model (\ref{eq:toymodelcov}).
  Coefficients $\left\{ \lambda^{(k)}_j \right\}_j$ of the generalized $\chi^2$ distribution in formula (\ref{eq:Q1})
  are eigenvalues of $\Sigma^{1/2} A_k \Sigma^{1/2}$, where the matrix $A_k$ is given by  (\ref{eq:Amatrix}).
  
  To illustrate obtained results for the FBM with noise we now compare the characteristic function of the estimator $\hat{r}(k)$ written as a quadratic form (\ref{eq:Q1def}) with the characteristic function of the generalized $\chi^2$ distribution (\ref{eq:Q1}). 
  The characteristic function of $\chi^2$ distribution with $k$ degrees of freedom is given by
  \begin{align*}
  \varphi_X(t) = \mathbb{E} \exp\{\mathrm{i} t X\} = (1- 2 \mathrm{i} t)^{-k/2},
  \end{align*}
  for $k>0$. The characteristic function of the quadratic form (\ref{eq:Q1}) is given by
  \begin{align*}
  \varphi_{Q}(t) = \prod_{j=1}^{N} \varphi_{\lambda^{(k)}_j U_j^2}(t) = \prod_{j=1}^{N} \varphi_{U_j^2}\left(\lambda^{(k)}_j t\right) = \prod_{j=1}^{N} \left(1-2 \mathrm{i} t \lambda^{(k)}_j\right)^{-1/2}.  
  \end{align*}
  
\begin{figure}[ht]
    \centering
      \includegraphics[width=0.6\textwidth]{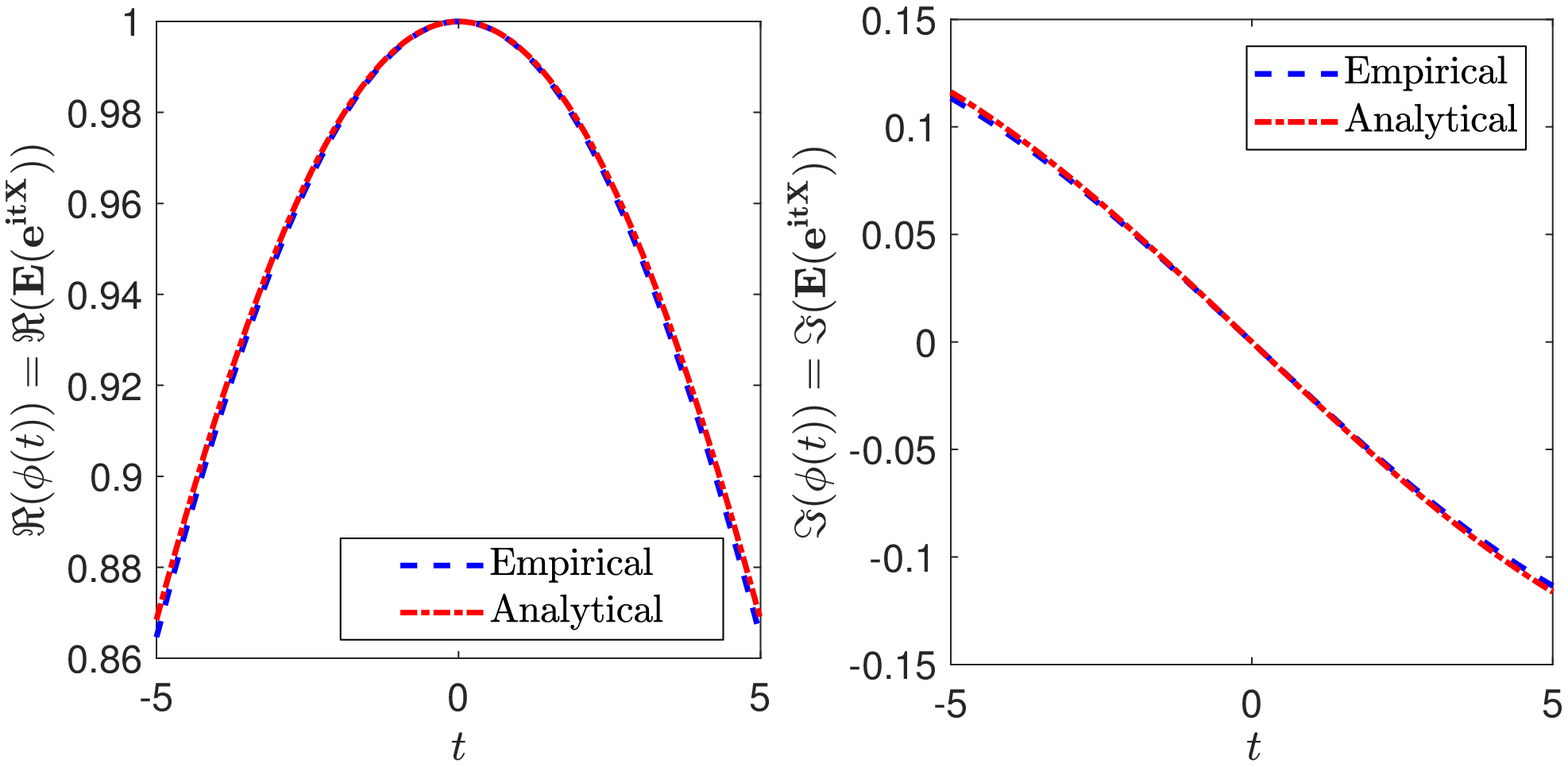}\\
      \includegraphics[width=0.6\textwidth]{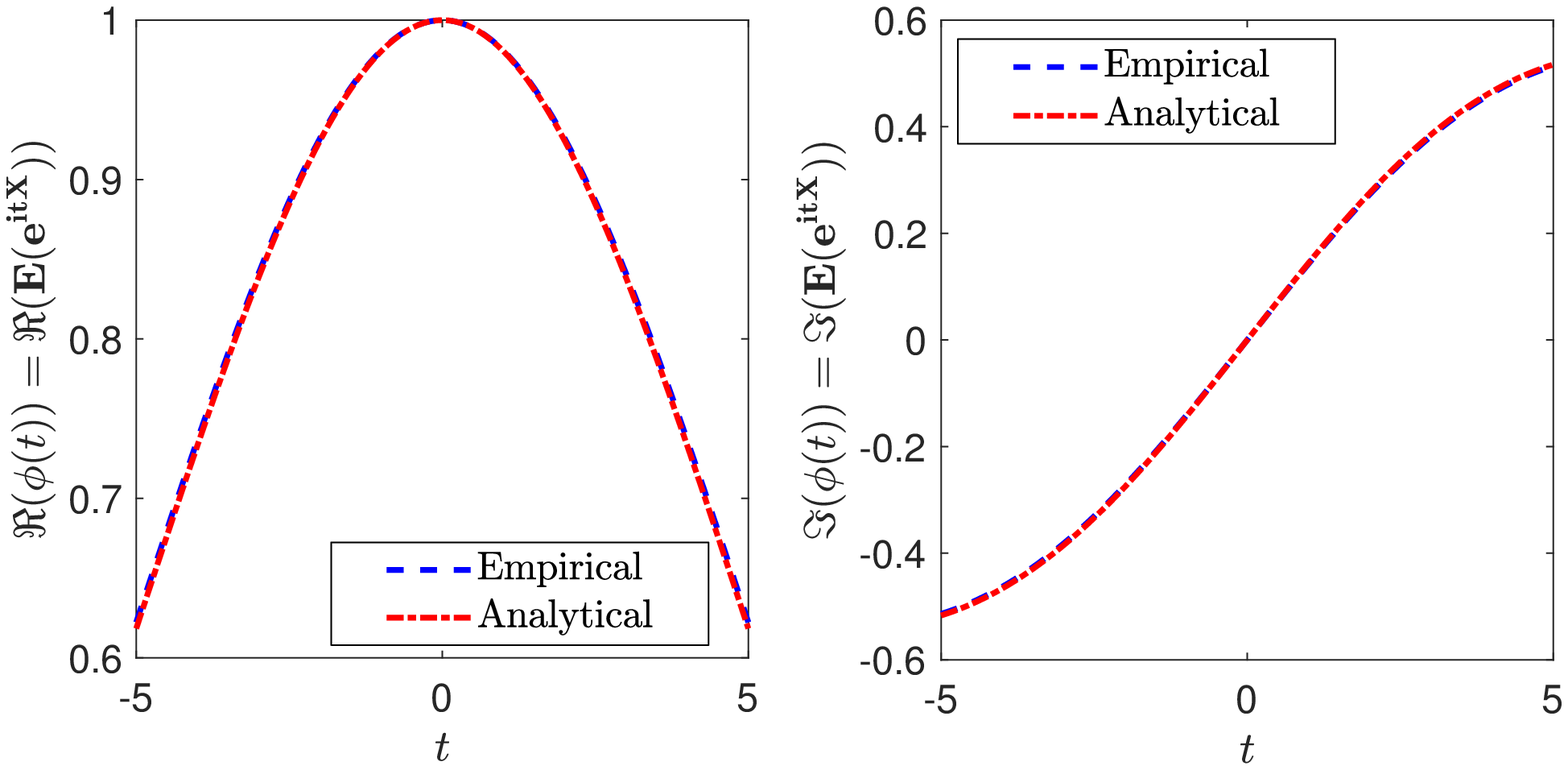}
    \caption{Comparison of the imaginary (left panel) and real (right panel) parts  of the empirical (dashed blue line) and analytical   (dash-dotted red line) characteristic functions of the estimator $\hat{r}(3)$ for the FBM and noise. The analytical distribution  of the estimator is given
    in terms of the generalized $\chi^2$ distribution. The top panel corresponds to a subdiffusion case with $H=0.3$ and the bottom panel to the superdiffusion with $H=0.7$. In both cases the magnitude of the additive noise is $\sigma=0.2$. The empirical characteristic function is calculated on the basis of $n=10^5$ trajectories of length $N = 2^7=128$.}
	\label{fig:rnComparison}
  \end{figure}

  Figure \ref{fig:rnComparison} presents a comparison of the empirical and analytical characteristic functions of the estimator $\hat{r}(k)$ for $k=3$.  The analytical CDF is given by the generalized $\chi^2$ distribution (\ref{eq:Q1}). The left panel presents the real part of the appropriate functions and the right panel the imaginary part. The top panel is related to a subdiffusion case with $H=0.3$, while the bottom to the superdiffusion with $H=0.7$. The characteristic function of $\hat{r}(k)$ is calculated by means of Monte Carlo simulations for $n=10^5$ and the data length $N=2^7$. 

\section{Test based on the autocovariance function estimator}
\label{sec:test}
In this section we propose a test on the FBM with noise based on the ACVF. The test will be based on Theorem \ref{th:Q1}, which describes the distribution of the sample ACVF for a model with a given covariance matrix. 


Specifically, we assume that null hypothesis is $\mathcal{H}_0:$ FBM with noise with $H=H_0$ and $\sigma=\sigma_0$, against $\mathcal{H}_1:$ it is not FBM with $H= H_0$ and $\sigma = \sigma_0$. The test statistic is given by (\ref{eq:rk}). Theorem \ref{th:Q1} states that under $\mathcal{H}_0$, the test statistic has the generalized $\chi^2$ distribution given by~(\ref{eq:Q1}). Thus, the critical set of the test, at significance level $a$, is given by $[q_{a/2},q_{1-a/2}]^c$, where $q_{a/2}$ and  $q_{1-a/2}$ are $a/2$ and $1-a/2$ quantiles of the distribution, respectively. An important question is what $k$ should be chosen in the test statistic? The answer, which is justified in the next section, is $k=1$.

\subsection{Construction of critical surfaces for the test}
\label{section:results}

In this part we construct a critical surface for the estimator for different self-similarity and noise magnitude parameters. We note that the following algorithm can be easily adapted to an arbitrary Gaussian process.


\begin{alg}{Algorithm to create the critical surface for the FBM with noise.}
  \begin{enumerate}
          \item {Choose:}
              \begin{itemize}
                  \item $N$ -- length of the trajectory;
                  \item $H$ -- Hurst index;
                  \item $\sigma^2$ -- variance of the measurement noise;                  \item $a$ -- significance level (i.e. we will be looking for quantiles of orders $a/{2}$ and $1-a/{2}$ of the estimator $\hat{r}(k)$).
              \end{itemize}
          \item Choose a specific lag $k$ for the ACVF. Then calculate the matrix $A_k$ given by (\ref{eq:Amatrix}). 
          \item Calculate the autocovariance matrix $\Sigma$ given by (\ref{eq:covMatrix}).
          \item Find the eigenvalues $\left\{\lambda_j^{(k)} \right\}_{j=1,\ldots, N}$ of the matrix $\Sigma^{1/2} A_k \Sigma^{1/2}$.
          \item Calculate quantiles of order $a/2$ and $1-a/2$ of the generalized $\chi^2$ distribution given by
              \begin{equation*}
                  \sum_{j=1}^N \lambda_j^{(k)} U_j^2,
              \end{equation*}
          where $U_j^2$'s are i.i.d. random variables with the $\chi^2$ distribution with one degree of freedom. They form top and bottom layers of the critical surface used for the testing purposes. We denote the critical surface by $ q(N, a,H, \sigma)$
  \end{enumerate}
  \end{alg}
  As a result, for a given trajectory length, we obtain a critical surface for the test as a function of the self-similarity parameter $H$ and magnitude of the noise $\sigma$.

\begin{figure}[ht!]
  \centering
    \includegraphics[width=0.6\textwidth]{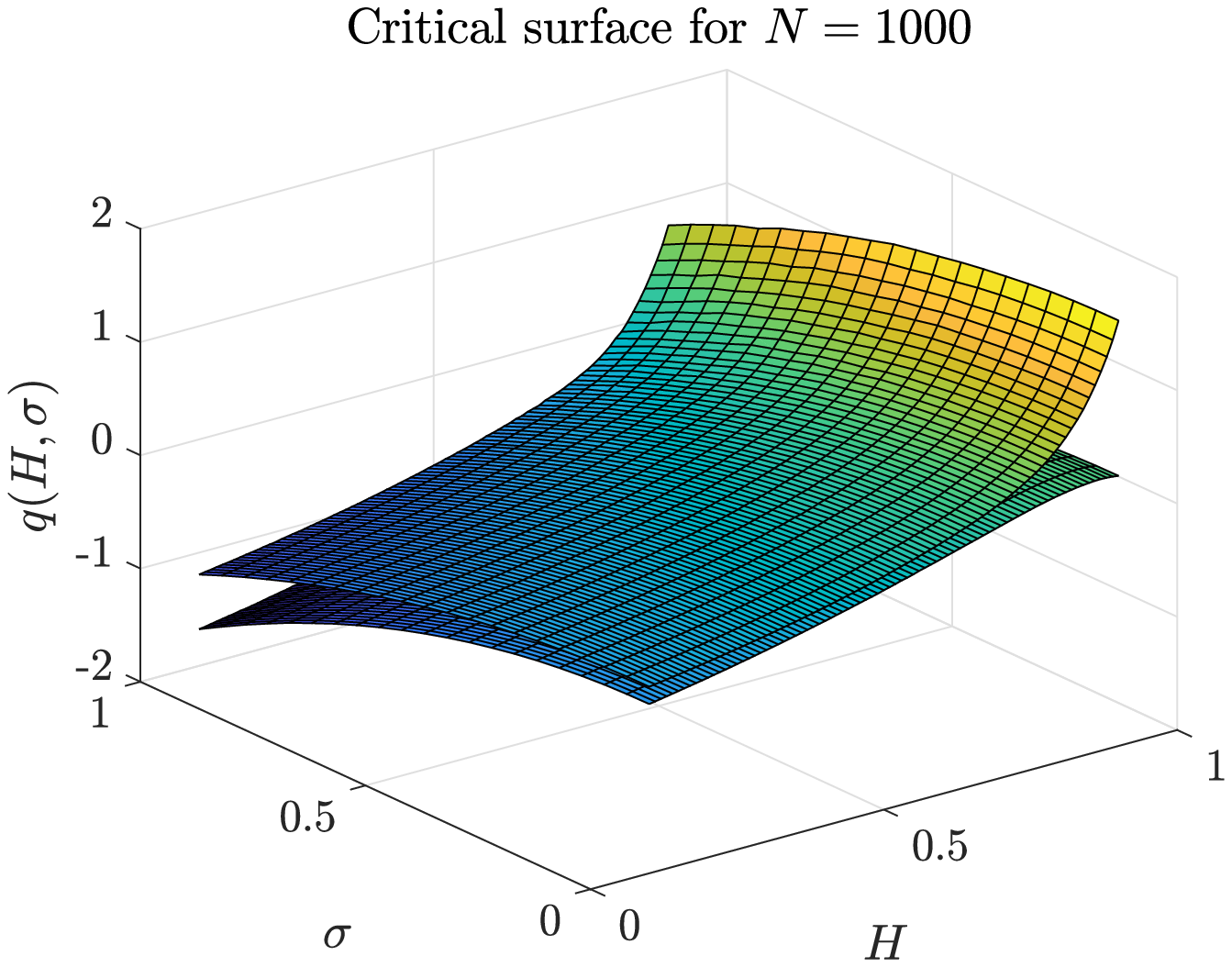}\\
    \includegraphics[width=0.47\textwidth]{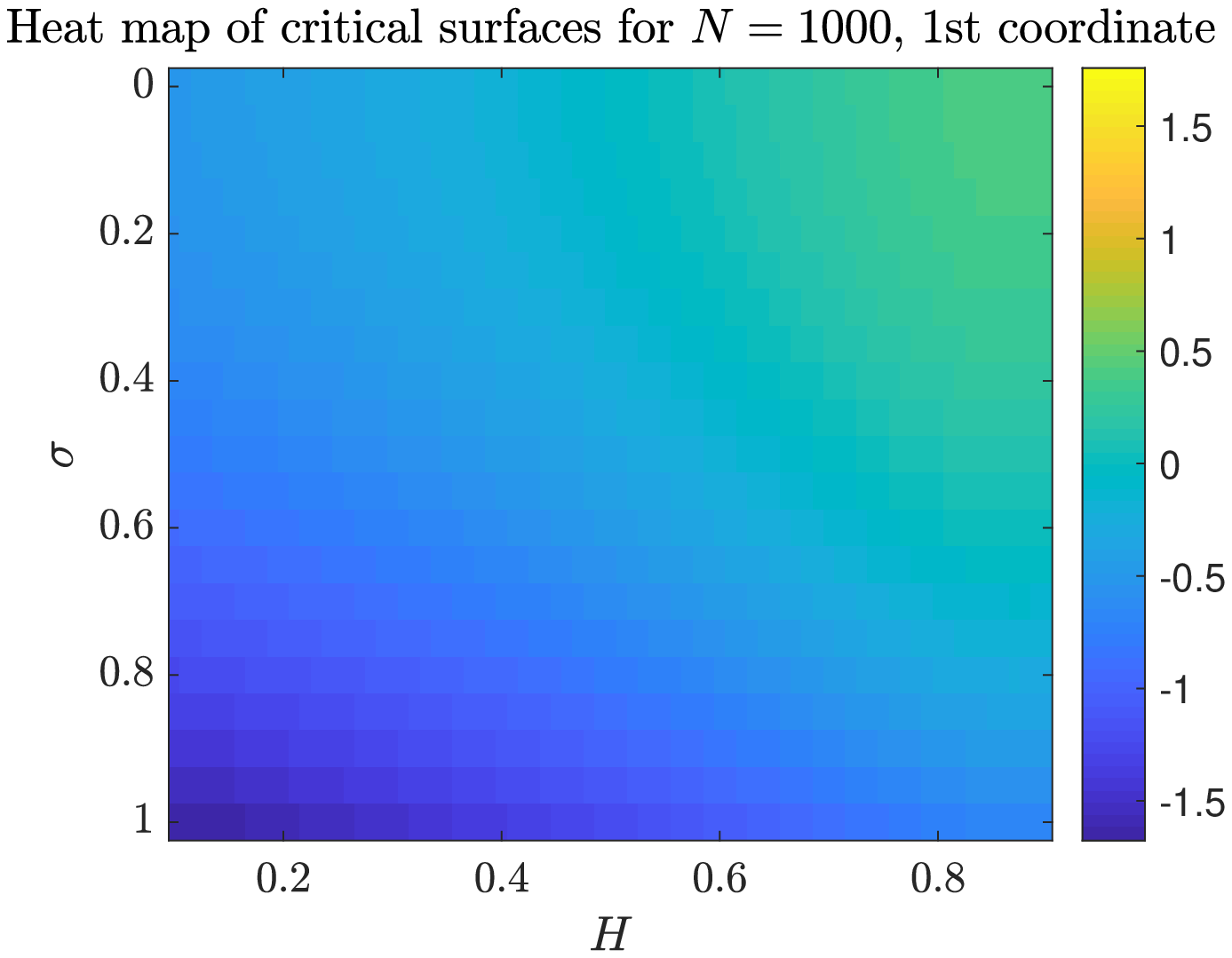}
    \includegraphics[width=0.47\textwidth]{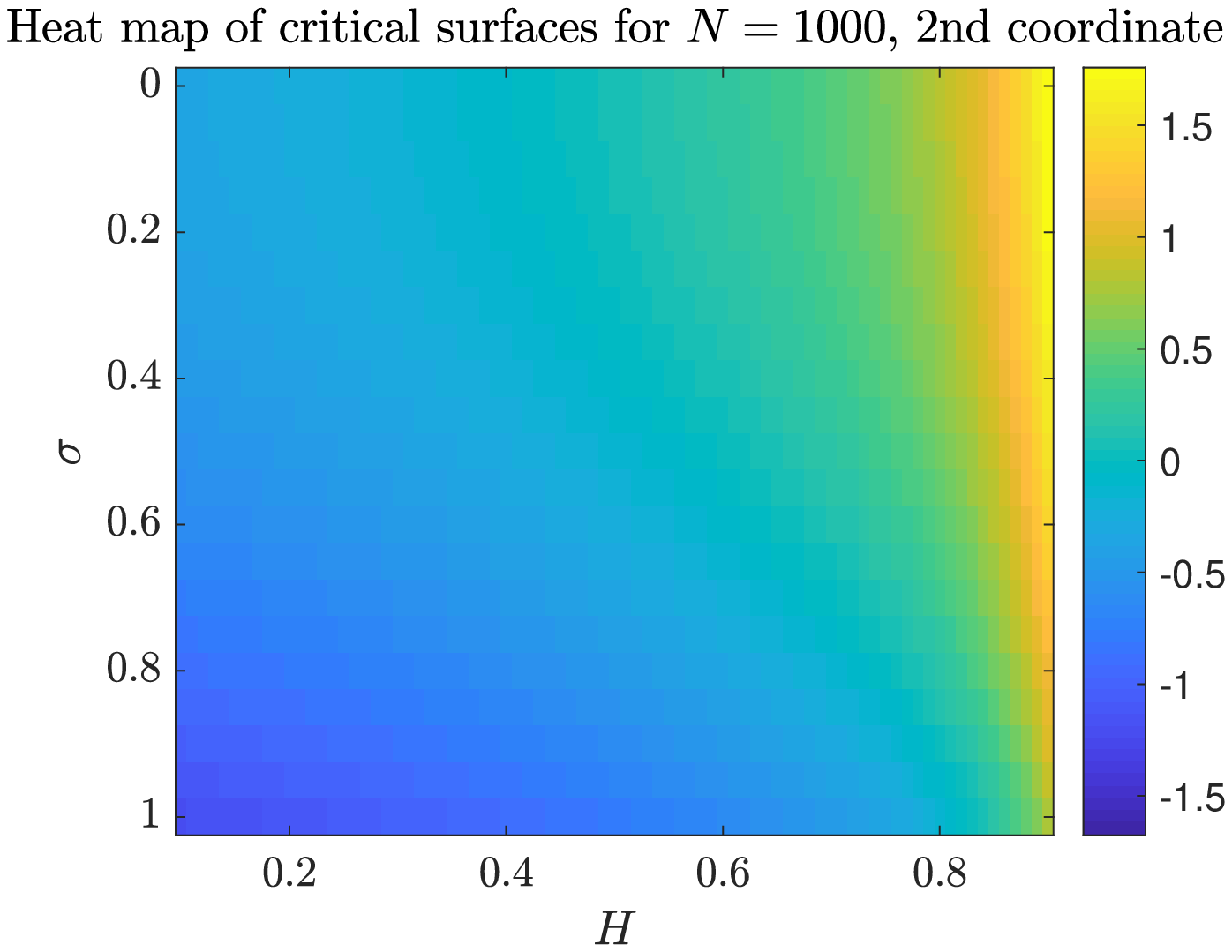}
    \caption{(Top panel) Critical surface $q(H, \sigma) = q(N=1000, a=0.05, H, \sigma)$ for parameters $H\in(0.1,0.9)$ and $\sigma\in[0,1]$. (Bottom panel) Heat maps for the lower (left panel) and upper (right panel) quantiles.   To create the surface we simulated $10, 000$ replications of the generalized $\chi^2$ distribution.}
    \label{fig:PEheatsurf1000}
\end{figure}

  In the top panel of Figure \ref{fig:PEheatsurf1000}  we present the critical surface $q(1000, 0.05,$ $H, \sigma)$ for parameters $H\in(0.1,0.9)$ and $\sigma\in(0,1)$. In the bottom panel of Figure~\ref{fig:PEheatsurf1000} we depict heat maps corresponding to this surface. We can notice that adding the additive noise described by the parameter $\sigma$ yields small changes, whereas big differences are caused by the changes in the $H$ parameter.



\begin{figure}[htp]
  \centering
    \includegraphics[width=0.6\textwidth]{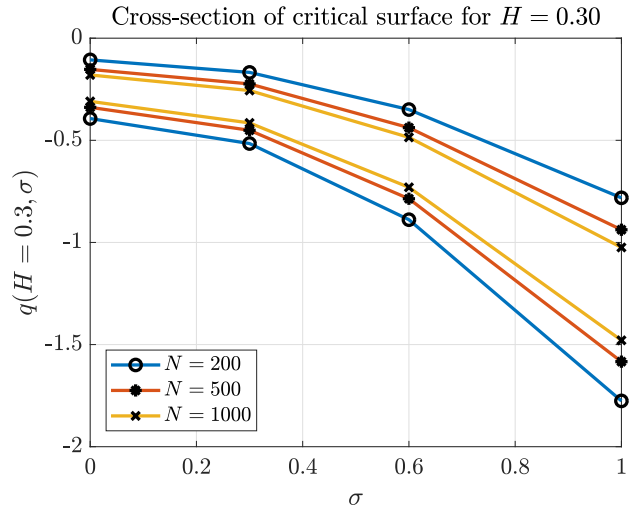} 
    \caption{Cross-section of the critical surface presenting quantile lines of the estimator (\ref{eq:rk}) for three data lengths $N$, $H=0.3$ and various $\sigma$'s. Blue lines correspond to length $N=200$, red to $N=500$ and yellow to $N=1000$. In each pair of lines of the  same colour the top line represents the quantile of order  $1-a/2=0.975$ and the bottom the quantile of order $a/2=0.025$.   To create the surfaces we simulated $10, 000$ replications of the generalized $\chi^2$ distribution.} 
    \label{fig:PEH0307}
\end{figure}

In Figure \ref{fig:PEH0307} we present quantiles $q_{0.05/2}$ and $q_{1-0.05/2}$ for three lengths of the trajectory: $N=200$ (blue), $N=500$ (red), and $N=1000$ (yellow), for the subdiffusive case $H=0.3$ as a function of the noise magnitude $\sigma$. For example, when the analysed data have length $N=200$ and we want to check if they come from FBM with noise with $H=0.3$ and $\sigma=0.3$, we should look at blue lines in Figure \ref{fig:PEH0307} at $\sigma=0.3$. We read the values -$0.51$ and -$0.16$. In such case, if the calculated value of the sample statistic given by $(\ref{eq:Q1def})$ lies between these numbers, then we do not have grounds for rejection that the data are described by the model with given parameters.

\begin{figure}[ht]
  \centering
    \includegraphics[width=0.6\textwidth]{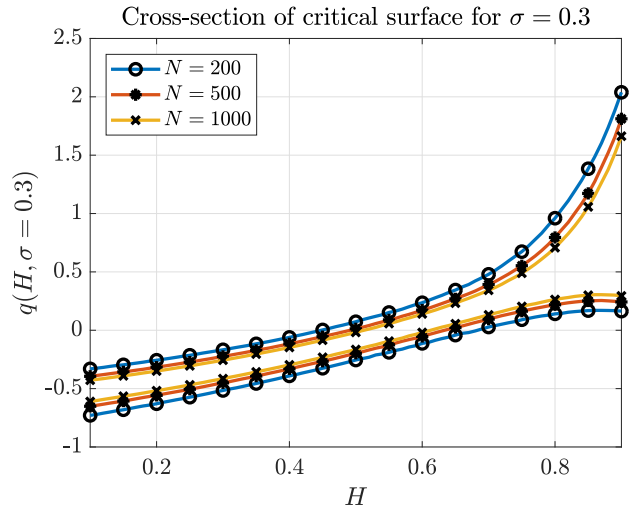}
    \caption{Cross-section of the critical surface presenting quantile lines of estimator (\ref{eq:rk}) for $\sigma=0.3$ and for three data lengths $N$.  Blue line correspond to  length $N=200$, red to $N=500$ and yellow to $N=1000$. In each pair of lines of the same colour the top line represents the quantile of order $1-\frac{a}{2}=0.975$, whereas the bottom the quantile of order $\frac{a}{2}=0.025$.}
    \label{fig:PEs02}
\end{figure}

In Figure \ref{fig:PEs02} we present functions $q_{0.05/2}$ and $q_{1-0.05/2}$ for four lengths of the trajectories: $N=200$ (blue), $N=500$ (red), and $N=1000$ (yellow), for the magnitude of the noise $\sigma=0.3$ as a function of the Hurst index $H$. For example, when the analysed data have length $N=200$ and $H=0.3$  we should look at blue lines in Figure \ref{fig:PEs02} and we can read the values -$0.39$ and -$0.11$. In such case, if the calculated value of the sample ACVF estimator given by $(\ref{eq:Q1def})$ lies between this numbers, then we do not have grounds to reject that the data are described by the model with given parameters.

We note that to estimate the Hurst exponent $H$ we can use a plethora of methods, e.g. Whittle estimator, detrended fluctuation analysis (DFA), rescaled range (R/S), DMA or MSD methods \cite{taqtev98,douetal03,beran2016long,sikora2018statistical,burgajsik11}.
For a method of estimation of both the self-similarity parameter and the magnitude of the measurement error for the considered model we refer the reader to Ref. \cite{burnecki2015estimating}.

\section{Power of the introduced test}
\label{sec:power}
In this section we present a Monte Carlo study done to show the power of the introduced test for models being FBM's with noise with different self-similarity parameters and noise magnitudes, and SBM's with different self-similarity parameters. We consider here three different null hypotheses corresponding to $H_0=0.3$, $0.5$ and $0.7$,
two different trajectory lengths ($N=200$ and $N=1000$) and for the FBM's two lags of the ACVF (1 and~2).
We carried out $n~=~10\,000$ replications and estimated the power for each model by dividing by $10\,000$ the number of times $H_0$ was rejected at $5\%$ significance level.

\begin{figure}[ht!]
  \centering
    \includegraphics[width=0.45\textwidth]{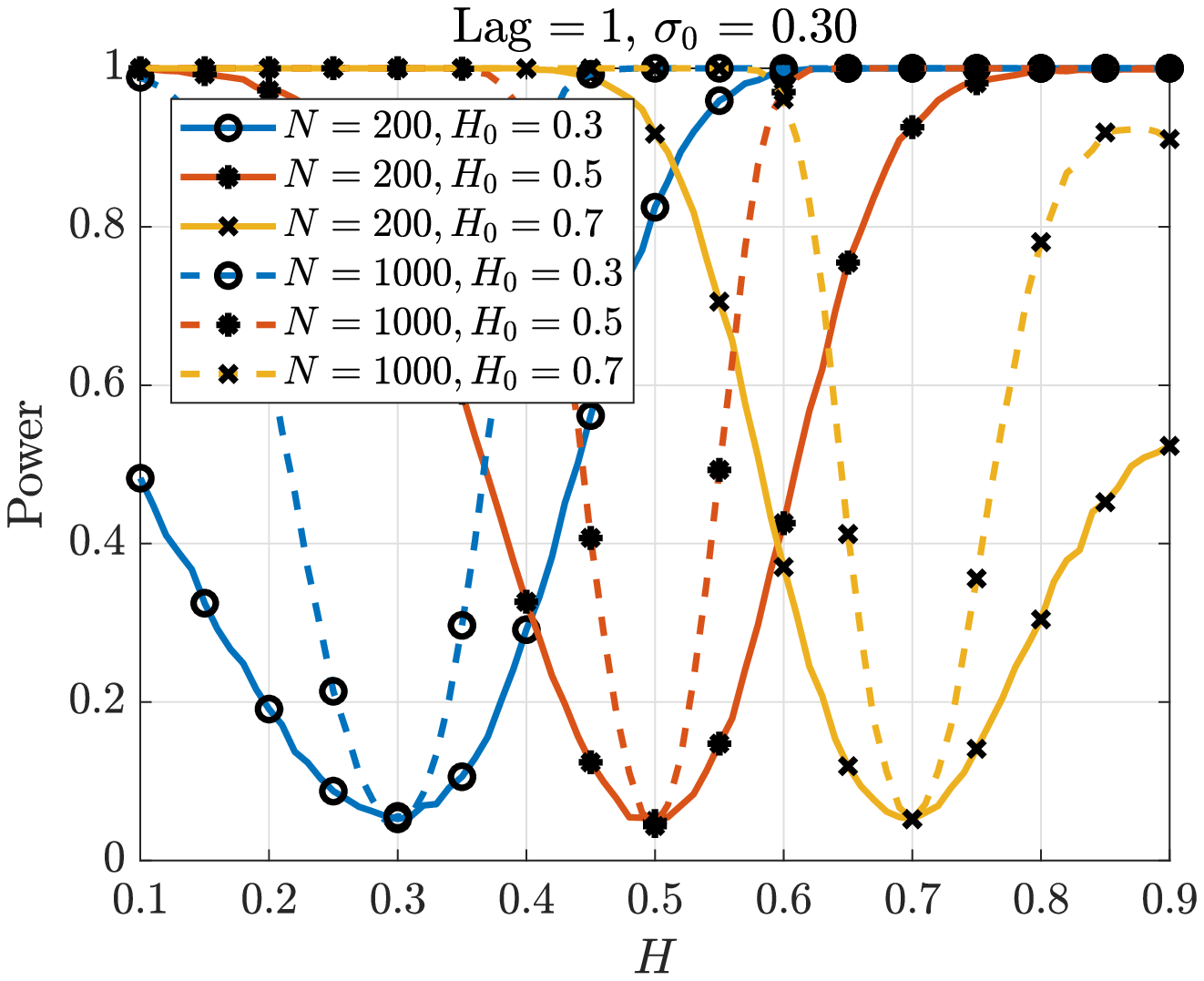}\\
    \includegraphics[width=0.45\textwidth]{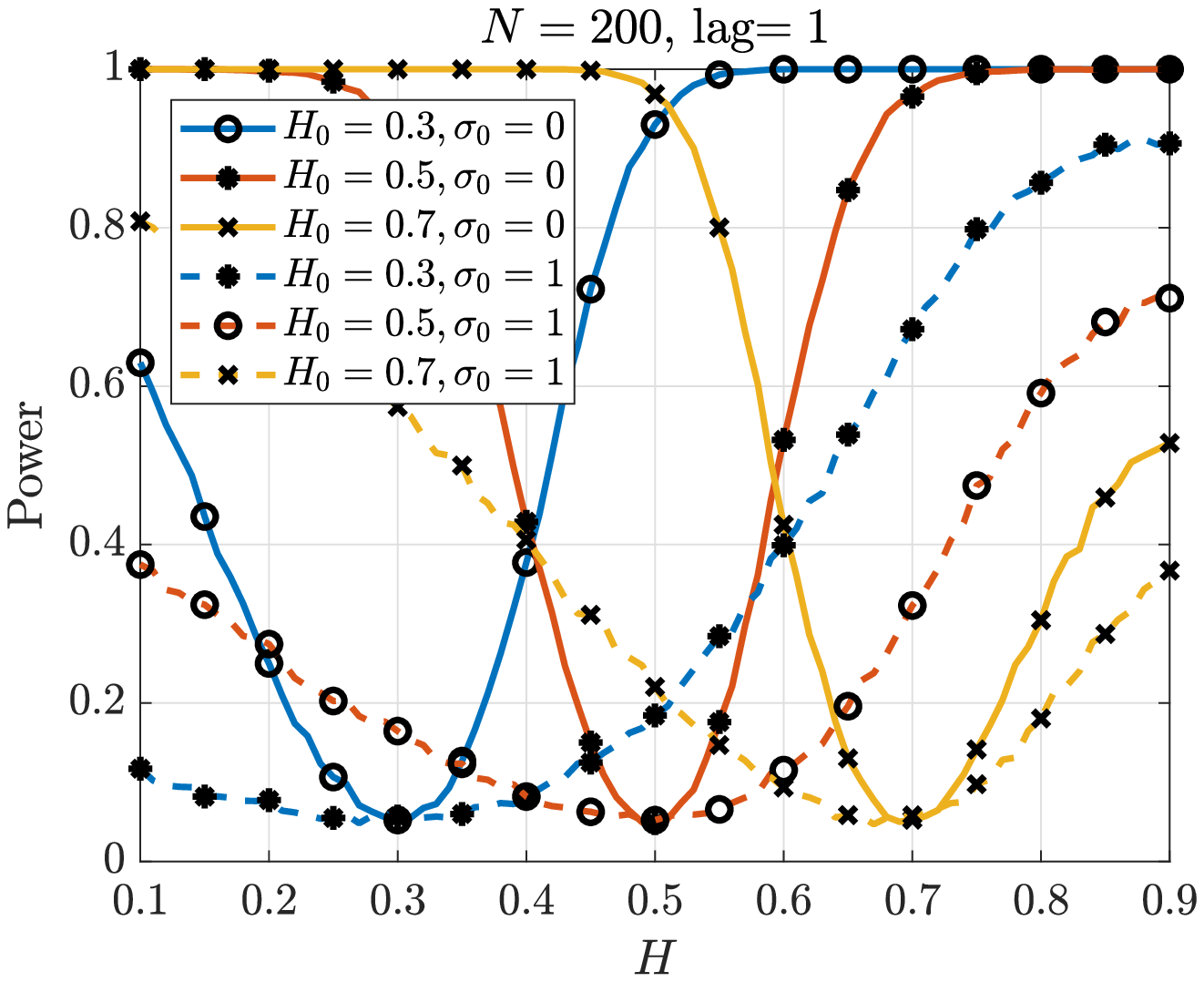}\\
    \includegraphics[width=0.45\textwidth]{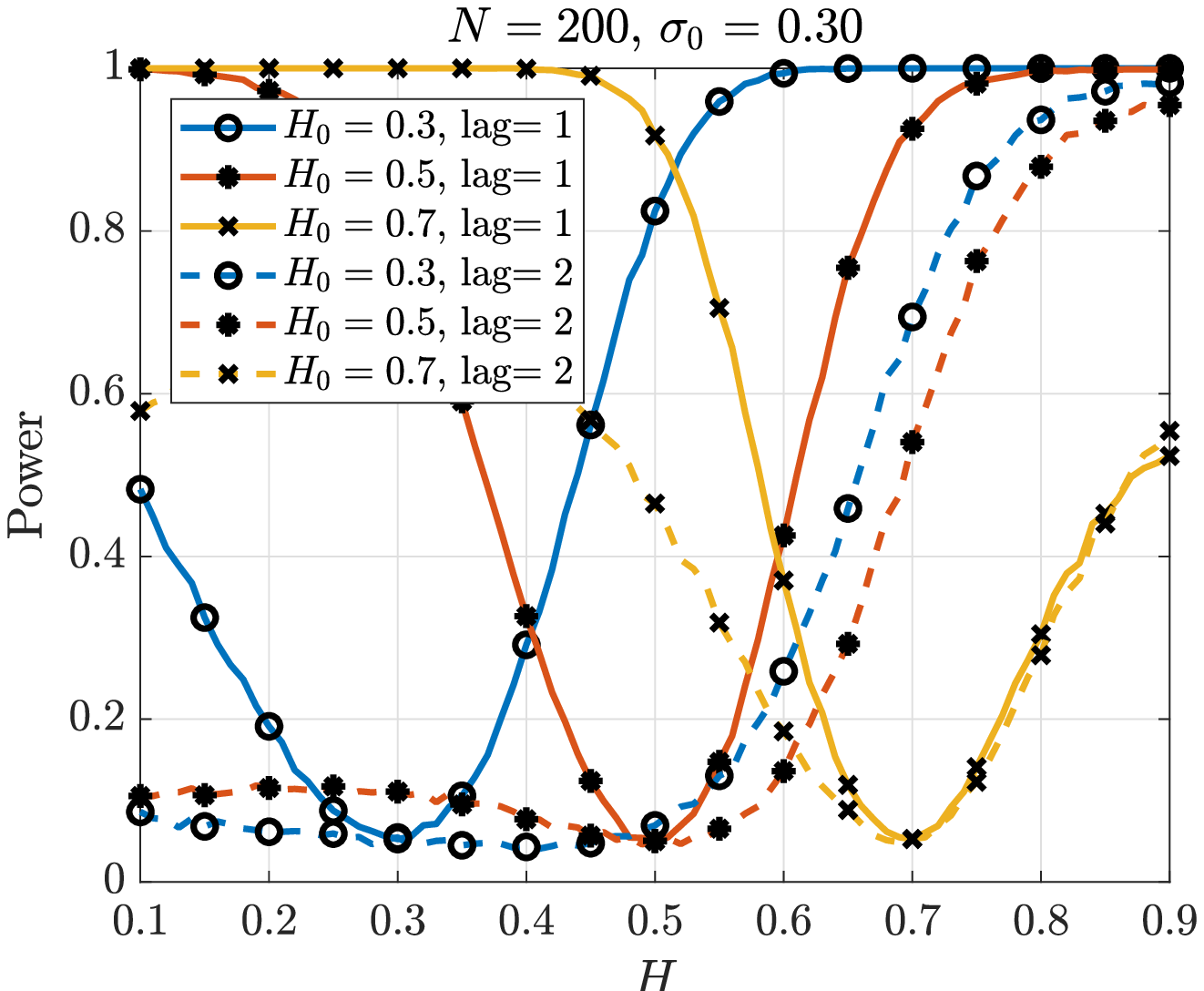}
    \caption{Test's power for different null hypotheses, noise magnitudes and data lengths, for simulated FBM's. Top panel: dependence on data length $N$ and Hurst index $H$, assuming $\sigma_0=0.3$. Middle panel: dependence on magnitude of the noise $\sigma$ and $H$, assuming $N=200$. Bottom panel: dependence on lag $\tau$ in the ACVF and $H$, assuming $\sigma_0=0.3$.}
    \label{fig:comparison}
\end{figure}

We start with an analysis of the power for models being FBM's with noise with different self-similarity parameters and noise magnitudes. We assume the same noise magnitude for both the null hypothesis and simulated model.
In the top panel of Figure \ref{fig:comparison} we illustrate an influence of length $N$ and Hurst index $H$ on the power of the test, assuming $\sigma_0=0.3$.
We can observe that the parameter $H_0$ seems not to have much impact on the power. We only note that the highest probabilities that the test rejects the false null hypothesis (the steepest parabola) are for $H_0=1/2$.
The middle panel presents a dependence of the power on both the magnitude of the noise $\sigma$ and $H$. We cab see that the power of the test decreases as $\sigma_0$ increases, the behaviour which is expected in a noisy environment since much noise makes it more difficult to distinguish between the models.
The bottom panel depicts a dependence on lag $\tau$ and $H$. 
We come to a conclusion that the power is much higher for lag $\tau=1$  than for $\tau=1$. We also checked the power of the test for lags up to 10 and the overall conclusion was the same, namely $\tau=1$ leads to the highest powers.

\begin{figure}[ht!]
  \centering
    \includegraphics[width=0.6\textwidth]{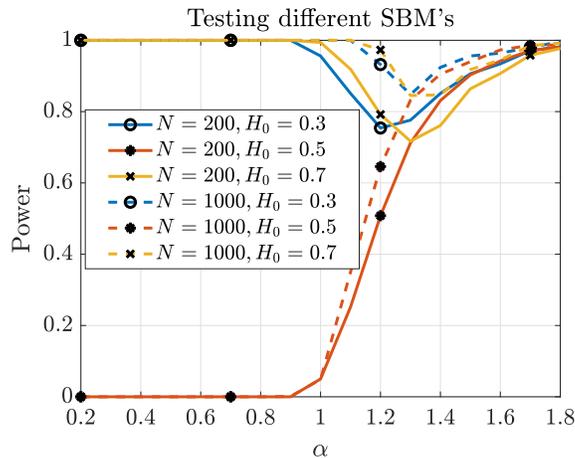}
    \caption{Test's power for different null hypotheses and data lengths, for simulated SBM's.}
    \label{fig:sbm-comparison}
\end{figure}

We now analyse the power for models being SBM's with different self-similarity parameters. For the null hypothesis we always take the noise magnitude equal to 0. In Figure \ref{fig:sbm-comparison} we illustrate the test's power. We note that the Hurst index of SBM $\alpha/2$ corresponds to $H$ of FBM so $\alpha=2H$. We can observe that for the null hypothesis with $H_0=0.3$ (corresponding to $\alpha=0.6$) power of the test slightly decreases for models with high $\alpha$'s, namely for $\alpha \in (1, 1.6)$ but is always very high. For the null hypothesis with $H_0=0.5$ (corresponding to $\alpha=1$) the test correctly distinguishes between SBM and FBM for high $\alpha$'s but for $\alpha<1$ the power is zero, thus it incorrectly does not reject the null hypothesis. For the case $H_0=0.7$ (corresponding to $\alpha=1.4$) the test's power is the lowest for $\alpha$'s corresponding to $2H$ but always exceeds 0.75.

\section{Conclusions}
\label{sec:con}

The FBM is a classical stochastic process to describe self-similar and long-range dependence phenomena. It has been applied to many different areas like such as telecommunication, economics, climatology and biology \cite{wiletal97,FALLAHGOUL201723,caretal07,metzler2014anomalous}. However, in many cases the recorded data are affected by a random noise which can be due, e.g., to the instrumentation error. We considered here a case of the white Gaussian noise added to the FBM. 
It is of great importance to be able to properly estimate the parameters of such process (see, e.g. \cite{burnecki2015estimating}) and to validate it.

In this paper we introduced a statistical test on the FBM with additive Gaussian noise based on the ACVF. We
derived a distribution of the test statistics which follows the generalised $\chi^2$ law. This allowed us to efficiently calculate critical surfaces, which were quantiles of the statistic as a function of the self-similarity parameter and magnitude of the noise. We presented an algorithm for a construction of critical surfaces for the FBM with noise
at a given significance level and different trajectory lengths.
We note that the procedure can be easily extended to an arbitrary Gaussian process with given parameters. 

We also note that two  tests for the pure FBM (without noise) were already proposed in the literature, see \cite{sikora2017mean, sikora2018statistical}. However, we emphasise simplicity of the introduced test (it is based directly on the ACVF which is known for many Gaussian processes) and that it accounts for the additive noise often present in the experimental data.

We checked the power of the introduced test by simulating alternatives being FBM's with different self-similarity parameters and noise magnitudes and SBM's with different self-similarity parameters.  We showed that the test can efficiently differentiate the studied model with given Hurst exponent and magnitude of the noise from other FBM's with noise. It can also distinguish between the FBM and SBM with the same self-similarity exponent.

\section*{Acknowledgements}
The authors acknowledge the support by NCN Maestro Grant No. 2012/ 06/A/ST1/00258. MB would like to additionally acknowledge the support of Wroc\l{}aw University of Science and Technology with the Research Grant No. 0402/0053/18.





\bibliographystyle{model1-num-names}
\bibliography{Bibliography.bib}







\end{document}